\documentclass[reqno]{amsart}
\usepackage{amsmath,amsthm,amssymb,hyperref,graphicx,mathrsfs,mathtools}  
\newtheorem{theorem}{Theorem}
\newtheorem{corollary}{Corollary}
\newtheorem{lemma}{Lemma}
\newtheorem{thmx}{Theorem}
\newtheorem{remark}{Remark}

\allowdisplaybreaks

\begin{document}
\title{ Certain Inequalities for the generalized polar derivative of a  polynomial}
\author{N. A. Rather}
\author{D. R. Bhat}
\author{I. Dar}
\address{$^{1,2,3}$Department of Mathematics, University of Kashmir, Srinagar-190006, India}
\email{dr.narather@gmail.com (Nisar Ahmad Rather), danishmath1904@gmail.com (Danish Rashid Bhat), ishfaq619@gmail.com (Ishfaq Ahmad Dar) }

\begin{abstract}
Recently Rather et al. \cite{NT} considered the generalized derivative and the generalized polar derivative and studied the relative position of zeros of generalized derivative and generalized polar derivative with respect to the zeros of polynomial. In this paper, we establish some inequalities that estimate the maximum modulus of generalized derivative and the generalized polar derivative of the polynomial $P(z)$, which is also the extension of recently known results.
\end{abstract}

\maketitle

 \section{\textbf{Introduction}}\normalfont
Let $P(z)$ be a polynomial of degree $ n $. A well known result, known as Bernstein's inequality \cite{S}, provides an upper bound for the maximum modulus of $ P'(z) $ on the unit circle $|z| = 1$ in terms of the maximum modulus of $ P(z) $, which is given by\\
\begin{align}
\displaystyle\max_{|z|=1} |P'(z)| \leq n \displaystyle\max_{|z|=1} |P(z)|.
\end{align}
Equality holds if and only if all the zeros of $ P(z) $ are at the origin.
For polynomials whose zeros lie entirely within the unit disk $|z| \leq 1$, Turán \cite{TUR} established a lower bound for the maximum modulus of $P'(z)$ on $|z| = 1 $. Specifically, if $P(z) $ is a polynomial of degree $ n $ with all zeros in $|z| \leq 1 $, then \\
\begin{align}
\displaystyle\max_{|z|=1} |P^\prime(z)| \geq \frac{n}{2} \displaystyle\max_{|z|=1} |P(z)|.
\end{align}
This inequality is sharp, and equality holds for $P(z) = \alpha z^n + \beta $, where $|\alpha| = |\beta| \neq 0 $.\\
\indent The estimation of the lower bounds of $Re\Big(\frac{zP^\prime(z)}{P(z)}\Big) $ on $|z|=1$, Dubinin \cite{V} proved the following.
\begin{thmx}\label{1.1}
If $P(z)=\displaystyle\sum_{j=0}^{n}a_jz^j$ has all its zeros in $|z|\leq1$, then for all $z$ on $|z| = 1$, for which $P(z)\neq0$,
\begin{align}
Re\Big(z\frac{P^{\prime}(z)}{P(z)}\Big) \geq \frac{n}{2}\Big\{1+\frac{1}{n}\Big(\frac{|a_n|-|a_0|}{|a_n|+|a_0|}\Big)\Big\}.
\end{align}
\end{thmx}
As an application of this, Dubinin \cite{V} obtained an interesting refinement of (2), by proving that if all the zeros of $P(z)$ lie in $|z|\leq1$, then
\begin{align}
\displaystyle\max_{|z|=1}|P^{\prime}(z)| \geq \frac{1}{2}\Big(n+\frac{|a_n|-|a_0|}{|a_n|+|a_0|}\Big)\displaystyle\max_{|z|=1}|P(z)|.
\end{align}
\indent Malik \cite{MA} extended $(2)$ and proved that if $P(z)\in P_n$ has all its zeros in $|z| \leq k$, $k \leq 1$, then 
\begin{align}\label{e2}
\underset{|z|=1}{max}|P^\prime (z)|\geq \frac{n}{1+k} \underset{|z|=1}{max}|P(z)|.
\end{align}
Equality in \eqref{e2} holds for the polynomial $P(z)=(z+k)^n$.\\
\indent Rather et al. \cite{NI} generalized Theorem A and proved the following.
 \begin{thmx}\label{1.2}\normalfont
 If $P(z)=\displaystyle\sum_{j=0}^{n}a_jz^j$ has all its zeros in $|z| \leq k$, $k \leq 1$, then for all $z$ on $|z| = 1$, for which $P(z)\neq0$,
\begin{align}
Re\Big(z\frac{P^{\prime}(z)}{P(z)}\Big) \geq \frac{n}{1+k}\Big\{1+\frac{k}{n}\Big(\frac{k^n|a_n|-|a_0|}{k^n|a_n|+|a_0|}\Big)\Big\},
\end{align}
which in turn yields the following refinement of the inequality $(5)$ as well as a generalization of the inequality $(4)$
\end{thmx}
\begin{thmx}\label{1.2}\normalfont
If $P(z)=\displaystyle\sum_{j=0}^{n}a_jz^j$ has all its zeros in $|z| \leq k$, $k \leq 1$ then for $|z| = 1$,
\begin{align}
|P^{\prime}(z)| \geq \frac{n}{1+k}\Big\{1+\frac{k}{n}\Big(\frac{k^n|a_n|-|a_0|}{k^n|a_n|+|a_0|}\Big)\Big\}|P(z)|.
\end{align}
\end{thmx}
Let $D_{\alpha}[P](z)$ denote the polar differentiation \cite{MM} of a polynomial $P(z)$ of degree $n$ with respect to a complex number $\alpha$, then
$$ D_{\alpha}[P](z)= nP(z)+(\alpha-z)P^{\prime} (z). $$
 Note that the polynomial $D_{\alpha}[P](z)$ is of degree at most $n-1$ and it generalizes the ordinary derivative $P^{\prime}(z)$ of $P(z)$ in the sense that
\[\underset{\alpha\rightarrow\infty}{\lim}\frac{D_{\alpha}[P](z)}{\alpha}=P^{\prime}(z)\] uniformly with respect $z$ for $\left|  z\right| \leq R,R>0.$\\
\indent The Bernstein-type inequalities for the class of polynomials with  “ordinary derivative” replaced by “polar derivative” have attracted number of mathematicians. In this direction, Aziz \cite{AI} was the first to establish inequalities concerning the polar derivative of a polynomial in terms of the modulus of the polynomial on the unit disk.\\
\indent As an extension of the inequality $(1)$ to the polar derivative, Aziz \cite{AI} proved that\\
if $P(z)$ is a polynomial of degree n, then for every $\alpha \in \mathbb C$ with $|\alpha|\geq 1,$ we have 
$$ |D_\alpha [P](z)|\leq n|\alpha z^{n-1}|~ \underset{|z|=1}{max}|P(z)|\quad for\quad |z|\geq 1.$$
The result is best possible and equality in the above inequality holds for $P(z)= cz^n,\quad c\neq 0.$\\
Concerning the class of polynomials having all zeros in $|z|\leq k$, Aziz and Rather \cite{AAA} obtained several sharp results concerning the maximum modulus of $D_{\alpha}[P](z)$ on $|z|=1.$  Among other things, they  established the following extension of the inequality \eqref{e2}  to the polar derivative of a polynomial.\\
\begin{thmx}\label{ta}
 If all the zeros of P(z) lie in $\left|z\right|\leq k$, then for $\alpha \in \mathbb C$ with $\left|\alpha\right|\geq k \leq 1, $  
\begin{align}\label{e4}
\underset{|z|=1}{max}|D_\alpha[P](z)|\geq \frac{n}{1+k}(|\alpha|-k)\underset{|z|=1}{max}|P(z)|.
\end{align}
\end{thmx}
Rather et al. \cite{NI} refined this inequality by proving the following result
\begin{thmx}\label{ta}
If all the zeros of $P(z)=\displaystyle\sum_{j=0}^{n}a_jz^j$ lie $|z| \leq k, k \leq 1$, then for  $\alpha \in \mathbb C$ with $\left|\alpha\right|\geq k ,$  
\begin{align}
\displaystyle\max_{|z|=1}|D_{\alpha}[P](z)| \geq n\Big(\frac{|\alpha|-k}{1+k}\Big)\Big\{1+\frac{k}{n}\Big(\frac{k^n|a_n|-|a_0|}{k^n|a_n|+|a_0|}\Big)\Big\}\displaystyle\max_{|z|=1}|P(z)|.
\end{align}
\end{thmx}
In literature there exist several generalizations and refinements of these inequalities. For $P(z)=c(z-z_1)(z-z_2)\cdots(z-z_n)$ and $\mathbb R_+^n$ being the set of all n-tuples $\gamma := (\gamma_1, \gamma_2,\ldots, \gamma_n)$ of non-negative real numbers (not all zeros) with $\sum\limits_{j=1}^{n}\gamma_j=\wedge$ , Sz-Nagy \cite{SN} introduced generalized derivative of polynomial $P(z)$ defined by
$$P^\gamma(z):=P(z)\displaystyle\sum_{j=1}^{n}\frac{\gamma_j}{z-z_j}= \sum\limits_{j=1}^{n} \gamma_j P_k(z). $$
Notice that  for $\gamma=(1,1,1,\ldots,1),P^\gamma(z)=P^\prime(z).$ In view of this, we  call it generalized derivative of polynomial $P(z)$.
\\Next we define generalized polar derivative of $P(z)$ as
$$D_{\alpha}^{\gamma}[P](z):=\wedge P(z)+(\alpha-z)P^\gamma(z),$$
where $\wedge =\sum\limits_{j=1}^{n} \gamma_j.$
\\Noting that for $\gamma=(1,1,1,\ldots,1), \quad D_{\alpha}^{\gamma}[P](z)=D_\alpha P(z)$, we call it generalized polar derivative of $P(z)$.\\
\indent Rather et al. \cite{NL} extended Theorem D to the generalized polar derivative and proved the following result.
\begin{thmx}\label{ta}
If all the zeros of P(z) lie in $\left|z\right|\leq k$, $k\leq1$ then for $\alpha \in \mathbb C$ with $\left|\alpha\right|\geq k , $  
\begin{align*}
\underset{|z|=1}{max}|D_\alpha^\gamma[P](z)|\geq \frac{\wedge}{1+k}(|\alpha|-k)\underset{|z|=1}{max}|P(z)|.
\end{align*}
\end{thmx}
Recently Rather et al. \cite{NT} proved the following results
\begin{thmx}\label{ta}\normalfont
If $P(z)=\displaystyle\sum_{j=0}^{n}a_jz^j$ is a polynomial of degree $n$ having  all its zeros in $\left|z\right|\leq k$ where $k\leq 1,$ then 
\begin{align}\label{theq}
|P^\gamma(z)|\geq \frac{k}{1+k}\Big[\frac{\wedge}{k}+\gamma_{m}\Big(\frac{k^n|a_n|-|a_0|}{k^n|a_n|+|a_0|}\Big)\Big]|P(z)|,
\end{align}
where $\gamma_m=\min\{\gamma_1, \gamma_2, \cdots, \gamma_n \}$.
\end{thmx}
\begin{thmx}\label{ta}
If $P(z)=\displaystyle\sum_{j=0}^{n}a_jz^j$ is a polynomial of degree $n$ having  all its zeros in $\left|z\right|\leq k$ where $k\leq 1,$ then for $\alpha \in \mathbb{C}$ with $|\alpha|\geq1,$
\begin{align*}
\underset{|z|=1}{max}|D_\alpha^\gamma [P](z)|\geq \Big(\frac{|\alpha|-k}{1+k}\Big)\Big[\wedge+k\gamma_m\Big(\frac{k^n|a_n|-|a_0|}{k^n|a_n|+|a_0|}\Big)\Big]\underset {|z|=1}{ max}|P(z)|,
\end{align*}\\
where $\gamma_m=\min\{\gamma_1,\gamma_2,\cdots,\gamma_n\}$ and the result is best possible for polynomial $P(z)=(z+k)^n$
\end{thmx}
For $k=1$, Theorem H reduces to the following result.
\begin{thmx}\label{ta}
If $P(z)=\displaystyle\sum_{j=0}^{n}a_jz^j$ is a polynomial of degree $n$ having  all its zeros in $\left|z\right|\leq 1$, then for $\alpha \in \mathbb{C}$ with $|\alpha|\geq1,$
\begin{align*}
\underset{|z|=1}{max}|D_\alpha^\gamma [P](z)|\geq \Big(\frac{|\alpha|-1}{2}\Big)\Big[\wedge+\gamma_m\Big(\frac{|a_n|-|a_0|}{|a_n|+|a_0|}\Big)\Big]\underset {|z|=1}{ max}|P(z)|,
\end{align*}\\
where $\gamma_m=\min\{\gamma_1,\gamma_2,\cdots,\gamma_n\}$ and the result is best possible for polynomial $P(z)=(z+1)^n$.
\end{thmx}

\section{\textbf{Main Results}}
In this section,  we extend Theorem G and Theorem H to the class of  generalized polar derivatives of the polynomial whose zeros lie in $|z|$$ \leq k$, $k \geq 1$. We begin by proving the following results.
\begin{theorem}\label{th1} \normalfont
If $P(z)=\displaystyle\sum_{j=0}^{n}a_jz^j$ is a polynomial of degree $n$ having  all its zeros in $\left|z\right|\leq k$ where $k\geq 1,$ then 
\begin{align}\label{theq}
\underset{|z|=1}{max}|P^\gamma(z)|\geq \frac{\wedge}{1+k^n}\Big[1+\frac{\gamma_{m}}{\wedge}\Big(\frac{k^n|a_n|-|a_0|}{k^n|a_n|+|a_0|}\Big)\Big]\underset{|z|=1}{max}|P(z)|,
\end{align}
where $\gamma_m=\min\{\gamma_1, \gamma_2, \cdots, \gamma_n \}$. The result is best possible for polynomial $P(z)=(z+k)^n$
\end{theorem}
If we take $\gamma=(1, 1, \cdots, 1)$ in Theorem $1$, we get the following result.
\begin{corollary}
If $P(z)=\displaystyle\sum_{j=0}^{n}a_jz^j$ is a polynomial of degree $n$ having  all its zeros in $\left|z\right|\leq k$ where $k\geq 1,$ then 
\begin{align}\label{theq}
\underset{|z|=1}{max}|P^\prime(z)|\geq \frac{n}{1+k^n}\Big[1+\frac{1}{n}\Big(\frac{k^n|a_n|-|a_0|}{k^n|a_n|+|a_0|}\Big)\Big]\underset{|z|=1}{max}|P(z)|.
\end{align}
\end{corollary}
\begin{theorem}\label{th2} 
If $P(z)=\displaystyle\sum_{j=0}^{n}a_jz^j$ is a polynomial of degree $n$ having  all its zeros in $\left|z\right|\leq k$ where $k\geq 1,$ then
\begin{align}\label{s}
\underset{|z|=1}{max}|D_\alpha^\gamma [P](z)|\geq \Big(\frac{|\alpha|-k}{1+k^n}\Big)\Big[\wedge+\gamma_m\Big(\frac{k^n|a_n|-|a_0|}{k^n|a_n|+|a_0|}\Big)\Big]\underset {|z|=1}{ max}|P(z)|,
\end{align}
where $\gamma_m=\min\{\gamma_1, \gamma_2, \cdots, \gamma_n \}$.
\end{theorem}
\begin{remark}\normalfont
If we divide inequality (\ref{s}) by $|\alpha|$ and then letting $|\alpha|\rightarrow \infty$, we get the inequality\\
\begin{align}
\underset{|z|=1}{max}|P^\gamma(z)|\geq \frac{\wedge}{1+k^n}\Big[1+\frac{\gamma_{m}}{\wedge}\Big(\frac{k^n|a_n|-|a_0|}{k^n|a_n|+|a_0|}\Big)\Big]\underset{|z|=1}{max}|P(z)|.
\end{align}
\end{remark}
If we take $\gamma=(1, 1, \cdots, 1)$ in Theorem $2$, we get the following result.
\begin{corollary}
If $P(z)=\displaystyle\sum_{j=0}^{n}a_jz^j$ is a polynomial of degree $n$ having  all its zeros in $\left|z\right|\leq k$ where $k\geq 1,$ then
\begin{align*}
\underset{|z|=1}{max}|D_\alpha [P](z)|\geq n\Big(\frac{|\alpha|-k}{1+k^n}\Big)\Big[1+\frac{1}{n}\Big(\frac{k^n|a_n|-|a_0|}{k^n|a_n|+|a_0|}\Big)\Big]\underset {|z|=1}{ max}|P(z)|.
\end{align*}
\end{corollary}
For the proof of the theorems we need following lemmas. The first is due to Rather et al. \cite{NI}.
\begin{lemma}\label{le5}
If  $0 \leq x_j \leq 1$, $j=1,2,\cdots,n,$ then
\begin{align*}
\displaystyle\sum_{j=1}^{n}\frac{1-x_j}{1+x_j}\geq \frac{1-\prod_{j=1}^{n}x_j}{1+\prod_{j=1}^{n}x_j} \qquad \forall  n \in N.
\end{align*}
\end{lemma}
By a simple consequence of Maximum Modulus theorem, we can get the following result.
\begin{lemma}\label{le5}
If $ P(z)$ is  a polynomial of degree n, then for $R\geq1$
\begin{align*}
\underset{|z|=R}{max}|P(z)|\leq R^n\underset{|z|=1}{max}|P(z)|.
\end{align*}
\end{lemma}
\begin{lemma}\label{le6}\cite{AAA}.
If $P(z)$ is a polynomial having all its zeros in the disk $|z|\leq k$ where $k\geq 1,$ then
\begin{align}\label{le61}
\underset{|z|=k}{max}|P(z)|\geq \frac{2k^n}{1+k^n} \underset{|z|=1}{max}|P(z)|.
\end{align}
\end{lemma}
\section{Proof of Theorems}
\begin{proof}[\bf{ Proof of Theorem \ref{th1}}] 
 Since $P(z)$ is a polynomial of degree $n$ having all its zeros in $|z | \leq k$, $k    \geq 1$.\\
  Let $G(z)=P(kz)$, then all the zeros of $G(z)$ lie in $|z| \leq 1$.\\
 Now
 \begin{align*}
G^{\gamma}(z)=G(z)\displaystyle\sum_{j=1}^{n}\frac{\gamma_j}{z-w_j} ,
 \end{align*}
where $w_j=\frac{z_j}{k}$ and $|w_j|\leq1$.\\
For $|z|=1,$
\begin{align*}
|zG^{\gamma}(z)|=& |G(z)|\Big|\displaystyle\sum_{j=1}^{n}\frac{z\gamma_j}{z-w_j}\Big|\\
               \geq& |G(z)|\displaystyle\sum_{j=1}^{n}\frac{\gamma_j}{1+|w_j|}\\
               =& \Big\{\displaystyle\sum_{j=1}^{n}\frac{\gamma_j}{2}\Big(1+\frac{1-|w_j|}{1+|w_j|}\Big)\Big\}|G(z)|\\
               =& \Big(\frac{\wedge}{2}+\frac{1}{2}\displaystyle\sum_{j=1}^{n}\frac{\gamma_j(1-|w_j|)}{1+|w_j|}\Big)|G(z)|\\
               \geq& \Big(\frac{\wedge}{2}+\frac{\gamma_m}{2}\displaystyle\sum_{j=1}^{n}\frac{1-|w_j|}{1+|w_j|}\Big)|G(z)|\\
               =& \Big(\frac{\wedge}{2}+\frac{\gamma_m}{2}\displaystyle\sum_{j=1}^{n}\frac{1-\frac{|z_j|}{k}}{1+\frac{|z_j|}{k}}\Big)|G(z)|\\
\end{align*}
Using Lemma 1, we get,
\begin{align*}               
  |zG^{\gamma}(z)|  \geq& \Big[\frac{\wedge}{2}+\frac{\gamma_m}{2}\Big(\frac{1-\displaystyle\prod_{j=1}^{n}\frac{|z_j|}{k}}{1+\displaystyle\prod_{j=1}^{n}\frac{|z_j|}{k}}\Big)\Big]|G(z)|\\
               =& \Big[\frac{\wedge}{2}+\frac{\gamma_m}{2}\Big(\frac{k^n|a_n|-|a_0|}{k^n|a_n|+|a_0|}\Big)\Big]|G(z)|.
\end{align*}
This mplies
\begin{align*}
\displaystyle\max_{|z|=1}\Big|G(z)\displaystyle\sum_{j=1}^{n}\frac{\gamma_j}{z-w_j}\Big| \geq \Big[\frac{\wedge}{2}+\frac{\gamma_m}{2}\Big(\frac{k^n|a_n|-|a_0|}{k^n|a_n|+|a_0|}\Big)\Big]\displaystyle\max_{|z|=1}|G(z)|\\
\end{align*}
Replacing $G(z)$ by $P(kz)$, we get,
\begin{align*}
 \displaystyle\max_{|z|=1}\Big|P(kz)\displaystyle\sum_{j=1}^{n}\frac{\gamma_j}{z-\frac{z_j}{k}}\Big| \geq \Big[\frac{\wedge}{2}+\frac{\gamma_m}{2}\Big(\frac{k^n|a_n|-|a_0|}{k^n|a_n|+|a_0|}\Big)\Big]\displaystyle\max_{|z|=1}|P(kz)|\\
 k\displaystyle\max_{|z|=k}\Big|P(z)\displaystyle\sum_{j=1}^{n}\frac{\gamma_j}{z-z_j}\Big| \geq \Big[\frac{\wedge}{2}+\frac{\gamma_m}{2}\Big(\frac{k^n|a_n|-|a_0|}{k^n|a_n|+|a_0|}\Big)\Big]\displaystyle\max_{|z|=k}|P(z)|
\end{align*}
Now from Lemma $3$, we have 
\begin{align*}
k\displaystyle\max_{|z|=k}\Big|P(z)\displaystyle\sum_{j=1}^{n}\frac{\gamma_j}{z-z_j}\Big| \geq \Big[\frac{\wedge}{2}+\frac{\gamma_m}{2}\Big(\frac{k^n|a_n|-|a_0|}{k^n|a_n|+|a_0|}\Big)\Big]\frac{2k^n}{1+k^n}\displaystyle\max_{|z|=1}|P(z)|.
\end{align*}
Thus we have,
\begin{align*}
 k\displaystyle\max_{|z|=k}\Big|P^{\gamma}(z)\Big| \geq \Big[\frac{\wedge}{2}+\frac{\gamma_m}{2}\Big(\frac{k^n|a_n|-|a_0|}{k^n|a_n|+|a_0|}\Big)\Big]\frac{2k^n}{1+k^n}\displaystyle\max_{|z|=1}|P(z)|.
\end{align*}
Using Lemma 2, we get,
\begin{align*}
 k^n\displaystyle\max_{|z|=1}\Big|P^{\gamma}(z)\Big| \geq \Big[\frac{\wedge}{2}+\frac{\gamma_m}{2}\Big(\frac{k^n|a_n|-|a_0|}{k^n|a_n|+|a_0|}\Big)\Big]\frac{2k^n}{1+k^n}\displaystyle\max_{|z|=1}|P(z)|.
\end{align*}
That is,
\begin{align*}
 \displaystyle\max_{|z|=1}\Big|P^{\gamma}(z)\Big| \geq \frac{\wedge}{1+k^n}\Big[1+\frac{\gamma_m}{\wedge}\Big(\frac{k^n|a_n|-|a_0|}{k^n|a_n|+|a_0|}\Big)\Big]\displaystyle\max_{|z|=1}|P(z)|\\
\end{align*}

This proves the theorem.

\end{proof}
\begin{proof}[\bf{ Proof of Theorem \ref{th2}}] 
Since the polynomial $P(z)$ has all its zeros in $|z| \leq k$ where $k \geq 1.$ \\ Therefore the polynomial $G(z)= P(kz)$ has all its zeros in $|z| \leq 1$.
\\Applying  Theorem I to the polynomial $G(z)$ and noting that $\big |{\frac{\alpha}{k}}\big| \geq 1,$ we get
\begin{align*}
\underset {|z|=1}{ max}\bigg|D_{\frac{\alpha}{k}}^\gamma [G](z)\bigg|\geq & \frac{\frac{|\alpha|}{k}-1}{2}\Big[\wedge+\gamma_m\Big(\frac{|k^na_n|-|a_0|}{|k^na_n|+|a_0|}\Big)\Big]\underset {|z|=1}{ max}|G(z)|\\
=& \frac{|\alpha|-k}{2k}\Big[\wedge+\gamma_m\Big(\frac{k^n|a_n|-|a_0|}{k^n|a_n|+|a_0|}\Big)\Big]\underset {|z|=1}{ max}|G(z)|.
\end{align*}
Replacing $G(z)$ by $P(kz)$, we get
\begin{align*}
\underset {|z|=1}{ max}\bigg|D_{\frac{\alpha}{k}}^\gamma [P](kz)\bigg| & \geq  \frac{|\alpha|-k}{2k}\Big[\wedge+\gamma_m\Big(\frac{k^n|a_n|-|a_0|}{k^n|a_n|+|a_0|}\Big)\Big]\underset {|z|=1}{ max}|P(kz)| \\ 
=&  \frac{|\alpha|-k}{2k}\Big[\wedge+\gamma_m\Big(\frac{k^n|a_n|-|a_0|}{k^n|a_n|+|a_0|}\Big)\Big]\underset {|z|=k}{ max}|P(z)|
\end{align*}
With the help of Lemma $3$, we have
\begin{align*}
\underset {|z|=1}{ max}\bigg|D_{\frac{\alpha}{k}}^\gamma [P](kz)\bigg| & \geq  \frac{|\alpha|-k}{2k}\Big[\wedge+\gamma_m\Big(\frac{k^n|a_n|-|a_0|}{k^n|a_n|+|a_0|}\Big)\Big]\underset {|z|=k}{ max}|P(z)| \\ 
& \geq  \frac{|\alpha|-k}{2k}\Big[\wedge+\gamma_m\Big(\frac{k^n|a_n|-|a_0|}{k^n|a_n|+|a_0|}\Big)\Big]\frac{2k^n}{1+k^n}\underset {|z|=1}{ max}|P(z)|.
\end{align*}
This gives
\begin{align}
\underset {|z|=1}{ max}\bigg|D_{\frac{\alpha}{k}}^\gamma [P](kz)\bigg| & \geq  \Big(\frac{|\alpha|-k}{1+k^n}\Big)\Big[\wedge+\gamma_m\Big(\frac{k^n|a_n|-|a_0|}{k^n|a_n|+|a_0|}\Big)\Big]k^{n-1}\underset {|z|=1}{ max}|P(z)|.
\end{align}
Also we have 
\begin{align*}
\underset {|z|=1}{ max}\bigg|D_{\frac{\alpha}{k}}^\gamma [P](kz)\bigg|& = \underset{|z|=1}{max} \bigg |\wedge P(kz)+\bigg(\frac{\alpha}{k}-z \bigg )P^\gamma(kz) \bigg |\\ &= \underset{|z|=1}{max} \bigg |\wedge P(kz)+\bigg(\frac{\alpha-kz}{k} \bigg )P(kz)\sum\limits_{j=1}^{n}\frac{\gamma_j}{z-\frac{z_j}{k} } \bigg |\\ &= \underset{|z|=1}{max} \bigg |\wedge P(kz)+\bigg(\frac{\alpha-kz}{k} \bigg )kP(kz)\sum\limits_{j=1}^{n}\frac{\gamma_j}{kz-z_j} \bigg |\\ &=\underset{|z|=1}{max} \bigg |\wedge P(kz)+(\alpha-kz )P(kz)\sum\limits_{j=1}^{n}\frac{\gamma_j}{kz-z_j} \bigg |\\ & =\underset{|z|=1}{max}  |F(kz)|\\ & =\underset{|z|=k}{max}  |F(z)|,
\end{align*}
where
$F(z)=\wedge P(z)+(\alpha-z )P(z)\sum\limits_{j=1}^{n}\frac{\gamma_j}{z-z_j}$
 is a polynomial of degree atmost $n-1$.
\\On using Lemma $2$, this gives
\begin{align*}
\underset {|z|=1}{ max}\bigg|D_{\frac{\alpha}{k}}^\gamma [P](kz)\bigg|& =\underset{|z|=k}{max}  |F(z)|\\& \leq k^{n-1}\underset{|z|=1}{max}  |F(z)|\\ & =k^{n-1}\underset{|z|=1}{max} \bigg |\wedge P(z)+(\alpha-z )P(z)\sum\limits_{j=1}^{n}\frac{\gamma_j}{z-z_j} \bigg |\\ &=k^{n-1}\underset{|z|=1}{max}|D_\alpha^\gamma [P](z)|.
\end{align*}
Hence, 
\begin{align}\label{theq2}
\underset {|z|=1}{ max}\bigg|D_{\frac{\alpha}{k}}^\gamma [P](kz)\bigg|\leq k^{n-1}\underset{|z|=1}{max}|D_\alpha^\gamma [P](z)|.
\end{align}
Combining inequalities $(14)$ and $(15)$, we get
\begin{align*}
k^{n-1}\underset{|z|=1}{max}|D_\alpha^\gamma [P](z)|\geq \Big(\frac{|\alpha|-k}{1+k^n}\Big)\Big[\wedge+\gamma_m\Big(\frac{k^n|a_n|-|a_0|}{k^n|a_n|+|a_0|}\Big)\Big]k^{n-1}\underset {|z|=1}{ max}|P(z)|.
\end{align*}
This implies,
\begin{align*}
\underset{|z|=1}{max}|D_\alpha^\gamma [P](z)|\geq \Big(\frac{|\alpha|-k}{1+k^n}\Big)\Big[\wedge+\gamma_m\Big(\frac{k^n|a_n|-|a_0|}{k^n|a_n|+|a_0|}\Big)\Big]\underset {|z|=1}{ max}|P(z)|.
\end{align*}
This proves the result.
\end{proof}

\end{document}